\documentclass [11pt,a4paper]{article}

\usepackage{amsmath} 

\usepackage{graphics}
\usepackage{latexsym} 
\usepackage{amssymb} 
\usepackage{enumerate} 
\def \N {\mathbb N}
\def \Z {\mathbb Z}

\def \R {\mathbb R}

\def \s {{\sigma}}
\def \sigm {{u}}

\def \prob        {\ensuremath{\mathbf{P}}}

\def \var         {\ensuremath{\mathbf{Var}}}
\def \cov         {\ensuremath{\mathbf{Cov}}}

\def\pt{\partial_t} 
\def\px{\partial_x} 

\def\us{\mathbf{s}} 
\def\uzero{\mathbf{0}} 
\def\uone{\mathbf{1}} 
\def\vareps{\varepsilon}

\title {Hydrodynamic Equation for a Deposition Model}
\author {B\'alint T\'oth
\and 
Wendelin Werner
}
\date {Technical University Budapest and Universit\'e Paris-Sud}
\begin {document}

\setlength{\baselineskip}{1.23\baselineskip}

\maketitle
%%%%%%%%%%%%%%%%%%%%%%%%%%%%%%%%%%%%%%%%%%%%%%%%%%%%%%%%%%%%%%%%%

\begin {abstract}
We show that the two-component 
system of hyperbolic 
conservation laws $\partial_t \rho + \partial_x (\rho u) 
=0 = \partial_t u + \partial_x \rho$
appears naturally in the formally 
computed hydrodynamic limit of 
some randomly growing interface models,
and we study some properties of this system.
\end {abstract}
\vskip 3mm

{\bf Key Words:} Hyperbolic conservation laws, KPZ equation

{\bf MSC-class.:} 35L65, 82C41, 60K35 
%
%%%%%%%%%%%%%%%%%%%%%%%%%%%%%%%%%%%%%%%%%%%%%%%%%%%%%%%%%%%%%%%%%

\section{Introduction}

The macroscopic behaviour of physical systems can 
often be described in  terms of non-linear partial 
differential equations. In many cases, it had  been 
shown that functionals of microscopic models 
from statistical physics converge in the hydrodynamic limit 
towards certain solutions of these partial differential equations.
 
Studying the partial differential equation
(or the system of partial differential equations) 
can turn out to be a  very hard challenge in itself:
Appearance of
singularities in finite time, shocks etc.  
The so-called hyperbolic conservation laws have in particular 
received a lot of interest. Even in one space dimension, 
these PDEs proved to be extremely interesting and challenging
both mathematically and phenomenologically. 
These are partial differential  equations of the type
\begin{equation*} 
\partial_t u + \partial_x J(u) = 0
\end{equation*}
where $u= u(t,x)$ 
takes its value in $\R^n$ and $J$ is a non-linear function 
from $\R^n$ into $\R^n$. 
 
The best known and most investigated examples are the following.
(See e.g. \cite{hormander, serre, smoller} for a comprehensive 
introduction and survey of the subject.) 
\begin{enumerate}[(1)]
\item
Burgers' equation (with no viscosity):  $n=1$ and 
\begin{equation*}
\label{burgers}
\pt u+ \px(u^2/2) = 0.
\end{equation*}
\item
The isentropic gas dynamics equation in one space dimension: 
$n=2$, the components are the density field $\rho(t,x)$ and momentum
field $m(t,x)$
\begin{equation}
\label{igd}
\left\{ 
\begin{array}{l}
\pt \rho +\px m=0
\\[5pt]
\pt m + \px \big( m^2/\rho +p(\rho)\big)=0
\end{array}
\right.
\end{equation}
where $p(\rho)$ is the pressure, depending on density only.
\item
The so-called p-system, which is an alternative formulation of
the dynamics of  one-dimensional gas. The two components
are the velocity field $u(t,x)$ and the specific volume (=
inverse density) field $v(t,x)$:
\begin{equation}
\label{psystem}
\left\{ 
\begin{array}{l}
\pt v - \px u=0
\\[5pt]
\pt u + \px p(v)=0.
\end{array}
\right.
\end{equation}
Here $p(v)$ denotes the pressure, as a function of specific
volume. 
\item
The shallow water equation is another two component system: 
$h(t,x)$ denotes the height of the (shallow) layer of water, 
$u(t,x)$ is the velocity field:
\begin{equation}
\label{shallowwater}
\left\{ 
\begin{array}{l}
\pt h + \px(h u) =0
\\[5pt]
\pt u + \px \big(u^2/2 +h\big)=0.
\end{array}
\right.
\end{equation}
\end{enumerate}

Since Riemann, a considerable amount of knowledge and technology 
(more recently, for instance, entropy solutions, compensated 
compactness method) has been derived that give a better 
understanding of the physically relevant solutions  to these
equations.

In the present paper, we will be considering a
particular  two-component 
(i.e. $n=2$) system of hyperbolic conservation laws  
that arises in the context of surface growth 
(or more precisely growing interfaces, since the surface is one-dimensional).
 In other words, at each time
$t \ge 0$, one sees a landscape $
x \mapsto h(t, x)$ where $x \in \R$. 
The function $h$ is increasing in time. The rough phenomenological 
description of the phenomena we are interested in corresponds to the
case where the surface is growing in the normal direction to
its boundary, but there exists a `tension' that tends to 
keep the surface together, in the sense that it will
fill in holes quickly.
In the physics literature, a famous equation has been 
proposed by Kardar, Parisi and Zhang (the KPZ equation)
as a model for such situations, cf. \cite{kardarparisizhang}. 
It is (in the mathematical 
jargon) an ill-posed non-linear partial differential
equation with a stochastic term:
\begin{equation*}
\label{kpzequ}
\pt h = \Delta h - (\px h)^2 + W
\end{equation*}
where $W=W(t,x)$ denotes  a space-time white noise.
We do not want to give a review of the huge physics
literature on this equation, but we briefly stress two 
aspects. 
(See \cite{barabasistanley} for a
state-of-the-art survey of the physics literature on the subject
and an exhaustive list of references up to 1995.)
First, there exists to our knowledge no completely satisfactory  
(see however \cite {bertinigiacomin}) derivation of 
this equation from a microscopic model. 
Second, it is predicted that 
`the' solution to this equation has a special scaling 
behaviour at late times. More precisely, it is believed
that when $\alpha, t, x$ are very large, 
$h^{(\alpha)} (t,x) = 
\alpha^{-1/3} h( \alpha t , \alpha^{2/3} x)$ 
is also a solution to the KPZ equation.
The exponents $1/3$ and $2/3$ should be related 
to various conjectures and recent rigorous results 
concerning the fluctuations of highest eigenvalues of
random matrices,  of first passage percolation
paths, of longest increasing sequences etc etc.

One  way to define one-dimensional 
interfaces $h(t,x)$ in terms of
particle systems goes as follows:  
Start with a (finite or infinite) system of particles    
that evolve randomly in the potential $h(t,x)$ (or in 
some potential defined in terms of $h$)  
and that all contribute to increase the potential 
in the sense that $h(t ,x)$ corresponds to the joint
local time (i.e. cummulated occupation time density) 
of the particles at time $t$ and site $x$. In other
words, $h(t,x)$ increases locally at $x$ if there is 
a particle at $x$ and time $t$.
Note that this leads naturally to a two-component system 
in  the (formally computed) hydrodynamical limit: 
the first component is the 
density of particles, and the second component is the 
gradient of the profile of the potential.

In \cite{tothwerner}, we constructed  a continuous stochastic
process, corresponding on a heuristic level to 
the case where there is exactly (and only) one particle (its location at
time $t$ is denoted by $X_t$)
which is  driven by  
\begin{equation*}
\label{tsrm1}
dX_t = - \px h (t, X_t) dt
\end{equation*}
and $h(t,x)$ is the local time of $X$ at $x$ and time $t$, 
so that 
\begin{equation*}
\label{tsrm2}
\pt h(t, x) = \delta(X_t -x).
\end{equation*}
For details concerning the construction and primary properties of
this process and a rigorous version of 
these equations, see \cite{tothwerner}. Let us just emphasize 
a couple of features: The process $(X_t, t \ge 0)$ is a random 
process, even though the previous `differential equations' 
look very deterministic.
One reason is that (in the stationary regime), the function 
$x \mapsto h(t,x)$ is not regular;  in fact, it is a 
Brownian motion in the space variable (for fixed $t$). 
Second, $X_t$ is not a usual stochastic process (it is 
not solution of a stochastic differential equation
for instance), it 
has the $2/3$ scaling: $(\alpha^{-2/3} X_{\alpha t}, t\ge 0)$ has
the  same law as $(X_t, t \ge 0)$. In particular,
$(\alpha^{-1/3} h(\alpha^{2/3}x, \alpha t), t \ge 0, x \in \R)$
has the same law as $(h(x,t), t \ge 0,
x \in \R)$ so that $h$ has the same scaling property as  
the asymptotic  scaling conjectured for the 
KPZ equation.  

The process $(X_t, t \ge 0)$ can be viewed as the scaling limit
of a discrete negatively reinforced (i.e. self-repellent) 
random walk $(S_n , n \ge 0)$ on $\Z$ called the 
`true self-avoiding walk' in the physics literature. This is 
a nearest-neighbour 
walk on $\Z$  that decides at each step to jump 
to the left or to the right according to a probability 
depending on how many times it has visited the neighbouring
sites (or edges) before. Suppose for instance that after 
$n$ steps $S_n=x$ and that the discrete walk
 $(S_i)_{i \le n}$ has jumped already $l$ (resp. $r$) times on the 
edge immediately to the left (resp. to the right) of $x$.  
Then, $S_{n+1} = x+1$ with probability 
\begin{equation*}
\label{tsaw}
\prob \big( S_{n+1} = x+1 \big |
l, r, S_n=x \big) 
= \frac {e^{-\beta l}} {e^{-\beta l} + e^{-\beta r}}
\end{equation*}
where $\beta>0$ is some fixed constant.
In other words, the walk will 
prefer to go along the edge it has visited less often in the
past.
Note also that the probability in fact  
depends only on the
difference $l-r$ (which depends on all the past trajectory).
The distribution of the rescaled position of the random walker, 
$S_n/n^{2/3}$, converges to (a multiple of) the one-dimensional
marginal 
distribution of the continuous process $X_t$ described above,
\cite{toth}. 

It seems natural to consider the case where this 
one particle is replaced by many particles performing the same
kind of self-repelling walk on $\Z$, with a \emph{joint 
cumulated local time of all particles}.  Or, in the continuous
space-time setting: a continuously distributed cloud of particles 
(that all contribute to the same local time), which is the subject of
the present paper. As we shall see, this leads in the (formally
computed) hydrodynamic limit to the following system of 
hyperbolic conservation laws:
\begin{equation}
\left\{
\begin{array}{l}
\pt \rho+\px (\rho u) =0
\\[5pt]
\pt u + \px \rho = 0
\end{array}
\right.
\label{ourpde}
\end{equation}
where $\rho$ corresponds to the density of particles at 
$x$ and time $t$, and $u (x,t) = -\px h$ corresponds to 
the negative gradient of the interface. 
It seems, that although this system looks very natural,
it has not been considered in the literature. 
We should emphasize that in spite of some formal similarities
with the p-system (\ref{psystem}) and the shallow water equation
(\ref{shallowwater}), the system (\ref{ourpde}) shows very
different behaviour and describes a quite different phenomenon.   
We hope that its study may lead to improved understanding
of some aspects of `growing interfaces' in general.
In particular, this equation could shed some light on some
of the conjectured properties of the KPZ equation. 
\emph{The goal of the present paper is not to present a complete
treatment of this system of partial differential equation,
but rather to initiate it as an alternative approach to 1-d
domain growth and deposition phenomena}.

\section{The PDE: phenomenological derivation} 

We define a deposition model in the following terms.
The actual state of the system is described by two functions:
\begin{equation*}
\rho:\R_+\times\R\to\R_+
\quad\mbox{ and }\quad
h:\R_+\times\R\to \R.
\end{equation*}
$\rho(t,x)$ is the density of the population performing the
deposition, while $h(t,x)$ is the deposition height at time $t$
and space coordinate $x$.  
The rules governing the time evolution of the system are the following
\begin{enumerate}[(1)]
\item
The total population is conserved, so
that the continuity equation
\begin{equation*}
\pt\rho+\px(\rho u)=0
\end{equation*}
is valid, where $u(t,x)$ is the velocity field, to be specified
by the dynamical rules.
\item
The deposition rate is proportional to  the density of the
population, i.e.
\begin{equation}
\pt h=c_1\rho,
\label{depositionrate}
\end{equation}
where $c_1$ is a positive constant.
\item
The population is driven by a velocity field proportional to
the negative gradient of height
\begin{equation}
u=-c_2\px h,
\label{velocityfield}
\end{equation}
where $c_2$ is another positive constant. 
This rule corresponds to the self-repellence mechanism 
described in the introductory section. 
\end{enumerate}
 From (\ref{depositionrate}) and (\ref{velocityfield}) we 
readily get 
\begin{equation*}
\pt u+c_1c_2\px \rho=0
\end{equation*}
Without loss of generality, we can
choose $c_1c_2=1$ and get the two component system of hyperbolic
conservation laws
\begin{equation}
\left\{
\begin{array}{l}
\pt \rho+\px (\rho u) =0
\\[5pt]
\pt u + \px \rho = 0
\end{array}
\right.
\label{ourpde2}
\end{equation}
This system of PDEs with initial conditions
\begin{equation}
\rho(0,x)=\rho^{(0)}(x),
\qquad
u(0,x)=u^{(0)}(x)
\label{ouric}
\end{equation}
is the main object of the present paper. 
%We
%propose it as a  hydrodynamic description of deposition
%phenomena, possible alternative of the stochastic KPZ approach. 

As a first remark we mention here the scale invariance of 
(\ref{ourpde2}). Let $\nu\in \R$ be fixed. Given the functions
$(t,x)\mapsto \rho(t,x)$ and $(t,x)\mapsto u(t,x)$ and a positive
fixed number $\alpha$, define the rescaled functions
\begin{align*}
\rho^{(\alpha)}(t,x)
&
:=
\alpha^{2(1-\nu)}\rho(\alpha t,\alpha^{\nu}x),
\\
u^{(\alpha)}(t,x)
&
:=
\alpha^{1-\nu} u(\alpha t,\alpha^{\nu}x).
\end{align*}
One can easily check that if $(\rho,u)$ is solution of
(\ref{ourpde2}), then 
$(\rho^{(\alpha)}, u^{(\alpha)})$ is also a solution,
for any $\alpha>0$. The choice $\nu=1$ yields the hyperbolic
scale invariance valid for any hyperbolic conservation law. More
interesting is for our purposes the choice $\nu=2/3$. This is
the physically relevant scale invariance, since the density
changes covariantly under this scaling, i.e., the total mass 
$\int\rho^{(\alpha)}dx$ is unchanged. 

With this choice of $\nu$ the following scale invariance of the
deposition height follows: 
\begin{equation*}
h^{(\alpha)}(t,x)
:=
\alpha^{-1/3} h(\alpha t,\alpha^{2/3}x).
\end{equation*}
Recall that this is exactly the \emph{conjectured asymptotic scale
invariance  of the one-dimensional KPZ equation}.

\section{Bricklayers}

We define a system of interacting particles living on $\Z$,
with \emph{two conserved} quantities, whose hydrodynamic modes
are governed by a two-component system of hyperbolic
conservation laws  which, after taking another limit (low
density/late time), transforms into our system (\ref{ourpde2}).
The computations of the present section are somewhat
formal. Working out all technical details (e.g. proving 
uniqueness of the equilibrium Gibbs measures or technical details
of Yau's hydrodynamic limit) needs more effort.  
The present section serves as microscopic motivation of the 
PDE proposed above. 

\subsection{The particle system}

The Great Wall of China is being built by a brigade of
bricklayers. The wall consists of columns of unit-size bricks,
piled above the edges of the lattice $\Z$. The  height of the 
column piled above  the edge $(j,j+1)$
(i.e., number of bricks in this column) is $h_j$. In the 
dynamics of the system the discrete negative gradients
$z_j:=h_{j-1}-h_{j}\in\Z$ will be relevant. 
The bricklayers occupy
the sites of the lattice. At each site $j\in\Z$ there might be
an unlimited number $n_j\in\N$ of bricklayers. Bricklayers jump
to neighbouring sites and at each jump $j\to j\pm1$ a brick is
added to the respective column of bricks. 

In more technical 
terms: particles (= bricklayers) perform continuous time nearest
neighbour walk on the lattice $\Z$ and $h_j$ measures the
cumulated (discrete) local time on the lattice edge $(j,j+1)$.

\begin{figure}[ht]
\centering
%\vspace {-2cm}
\resizebox{0.95 \textwidth}{!} {\includegraphics {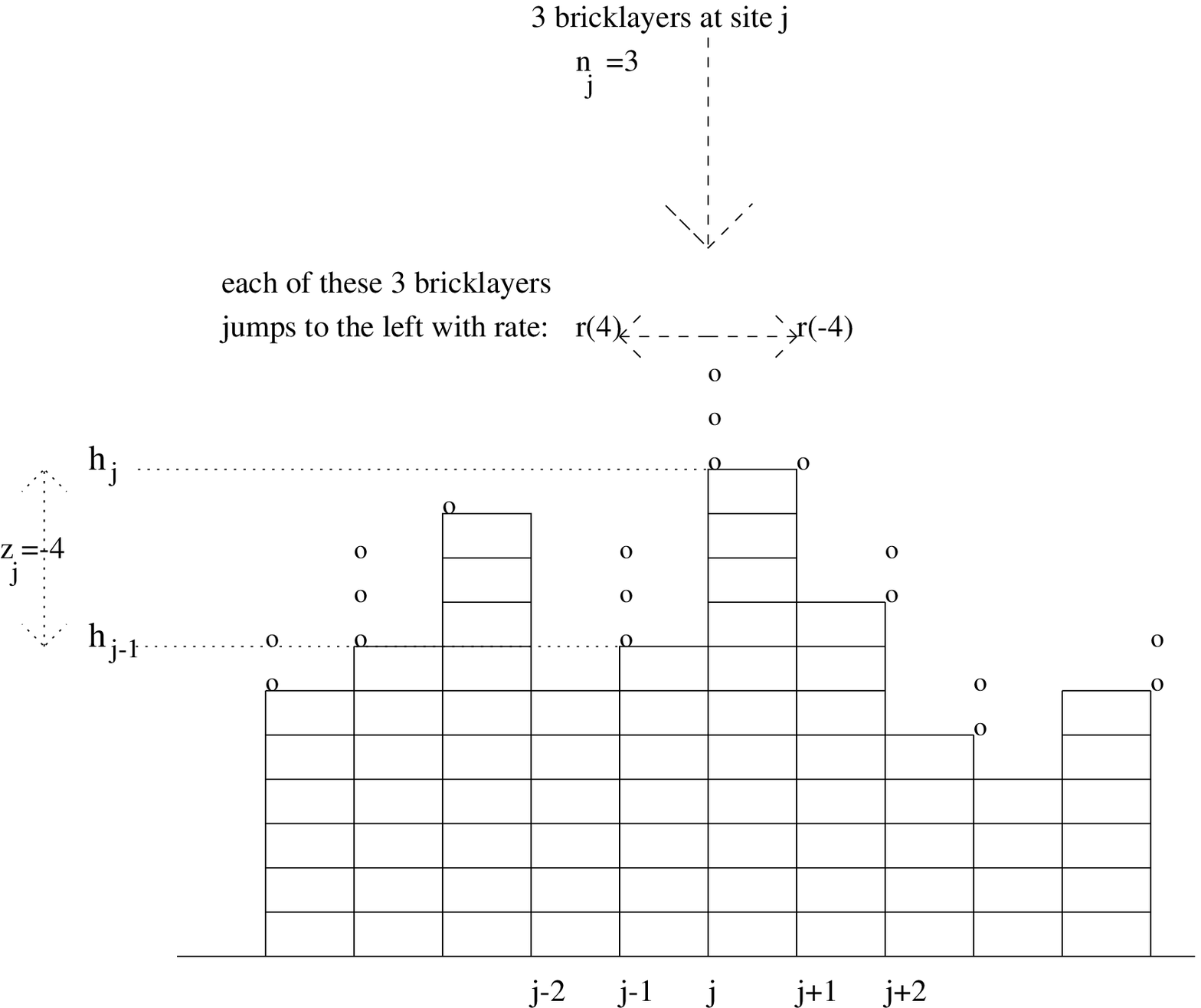}}
%\vspace {-1cm}
\caption{
The great wall being built}
\end{figure}

About the dynamics: the jump rates are chosen so that the
following conditions hold: 
\begin{enumerate}[(1)]
\item
the bricklayers' jumps are driven
by the local shape of the wall so that they try to reduce the
differences $z_j$ (i.e. to keep the height of the wall even),
\item
conditionally on the actual shape of the wall the
bricklayers jump independently. 
\end{enumerate}
% NEW 
This is done as follows. The instantaneous rate of jump from 
site $j$ to site $j\pm1$ (for each bricklayer sitting at site $j$)
is equal to $r(\pm z_j)$, where $r:\Z\to (0,\infty)$ is a fixed monotone
increasing function which defines the model. In order to be able
to compute  explicitely 
the stationary measures (see subsection 3.2) we impose
that $r(1-z) r(z)$ is a positive constant (this is
for instance the case if $r(z) = \exp (\beta z)$), and
multiplying time by a constant term, we can
in fact restrict ourselves
to the case where  
% END NEW
\begin{equation}
\label{rcond}
r(z)r(-z+1)=1, 
\quad
\text{ for all } 
z\in\Z.
\end{equation}
Thus, the following
changes of configuration may occur:
\begin{equation*}
(n_j,z_j),(n_{j+1},z_{j+1})
\to
(n_j-1,z_j-1),(n_{j+1}+1,z_{j+1}+1)
\end{equation*}
with rate $n_j r(z_j)$, and 
\begin{equation*}
(n_j,z_j),(n_{j-1},z_{j-1})
\to
(n_j-1,z_j+1),(n_{j-1}+1,z_{j-1}-1)
\end{equation*}
with rate $n_j r(-z_j)$.
% TO HERE
   
Clearly, $\sum_j n_j$ and $\sum_j z_j$ are formally conserved
quantities of the dynamics. It is also clear that besides these
globally conserved quantities the parity of $n_j+z_j$ is also
conserved on each lattice site $j\in\Z$. 

Now we give a more formal description of our interacting
particle system. For $s\in\{0,1\}$ let
\begin{equation*}
\big(\N\times\Z\big)_s:=
\{(n,z)\in\N\times\Z: n+z=s \mbox{ mod } 2\}.
\end{equation*}
Given the sequence $\us=(s_j)_{j\in\Z}\in\{0,1\}^{\Z}$ we define
the state space of our system as
\begin{equation*}
\Omega_{\us}:= 
\prod_{j\in\Z}\big(\N\times\Z\big)_{s_j}.
\end{equation*}
Elements of $\Omega_{\us}$ will be denoted by $\mathbf{\omega}$, i.e.
$\mathbf{\omega} = (\omega_j)_{j\in\Z}$ with 
$\omega_j = (n_j,z_j) \in \big(\N\times\Z\big)_{s_j}$. 
The (formal) infinitisimal generator of the Markov process
described  verbally in the first paragraph of this section, is: 
% FROM HERE
%
\begin{equation*}
Lf(\omega) =
\sum_{j\in\Z} n_j r(z_j)
\big(f(\Theta_{j+}\omega) - f(\omega)\big)
+
\sum_{j\in\Z} n_j r(-z_j)
\big(f(\Theta_{j-}\omega) - f(\omega)\big),
\end{equation*}
%
% TO HERE 
where the maps $\Theta_{j+}$ and $\Theta_{j-}$ act on the subsets
$\{\mathbf{\omega}\in\Omega_{\us}: n_j\ge1\}$ as
\begin{equation*}
\big(\Theta_{j+}\mathbf{\omega}\big)_i:=
\left\{
\begin{array}{lcl}
(n_i,z_i) 
&\quad \mbox{ if } \quad&
i\not=j, j+1
\\[5pt]
(n_i-1, z_i- 1) 
&\quad \mbox{ if } \quad&
i=j
\\[5pt]
(n_{i}+1, z_{i}+ 1)
&\quad \mbox{ if } \quad& 
i=j+1
\end{array}
\right.
\end{equation*}
respectively
\begin{equation*}
\big(\Theta_{j-}\mathbf{\omega}\big)_i:=
\left\{
\begin{array}{lcl}
(n_i,z_i) 
&\quad \mbox{ if } \quad&
i\not=j, j-1
\\[5pt]
(n_i-1, z_i+ 1) 
&\quad \mbox{ if } \quad&
i=j
\\[5pt]
(n_{i}+1, z_{i}- 1)
&\quad \mbox{ if } \quad&
i=j -1
\end{array}
\right.
\end{equation*}

\subsection{Equilibrium Gibbs measures}

% NEW 
For $k\ge0$ denote 
\begin{equation*}
R(z):=\prod_{k=1}^{|z|} r(k)
\end{equation*}
and 
\begin{equation*}
\theta^*:=\lim_{k\to\infty} r(k)\in(1,\infty].
\end{equation*}
Note that (\ref {rcond}) implies that for all $z \in \Z$, 
\begin {equation}
\label {rprop}
R(-z)=
R(z) = R(z-1) r(z) 
= R(z+1) r (-z).
\end {equation}
Fix the parameters $s\in\{0,1\}$, $\lambda>0$,
$\theta\in(1/\theta^*, \theta^*)$ and
define the probability measure $\mu_{s,\lambda,\theta}$ on
$\big(\N\times\Z\big)_s$ as follows:
\begin{equation*}
\mu_{s,\lambda,\theta}(n,z):=
\frac{1}{Z_s(\lambda,\theta)}
\frac{\lambda^n}{n!}\frac{\theta^z}{R(z)}, 
\end{equation*}
where 
\begin{equation*}
Z_s(\lambda,\theta):=
\sum_{(n,z)\in(\N\times\Z)_s}
\frac{\lambda^n}{n!}\frac{\theta^z}{R(z)}
\end{equation*}
is the normalizing factor (partition function). The measure
$\mu_{s,\lambda,\theta}$ is a product measure on $\N\times\Z$  
restricted to the subset $n+z=s\mbox{ mod }2$.
It is worth noting that  
\begin{equation}
\label{symmetries}
Z_s(\lambda,\theta)=
Z_s(\lambda,\theta^{-1})
\quad
\text{ and }
\quad
\mu_{s,\lambda,\theta}(n,z)
=
\mu_{s,\lambda,\theta^{-1}}(n,-z).
\end{equation}

For a fixed sequence $\us\in\{0,1\}^\Z$ and fixed parameters
$\lambda>0$, $\theta\in(1/\theta^*,\theta^*)$ 
we define on $\Omega_{\us}$ the
probability measure  
\begin{equation*}
\mu_{\us,\lambda,\theta}:=
\prod_{j\in\Z} \mu_{s_j, \lambda, \theta}.
\end{equation*}

By direct computations, one can check using (\ref {rprop})
that 
for any function $f$
that depends only on the value
of finitely many ${\mathbf {\omega}}_k$'s, for
any fixed $s_j$, $s_{j+1}$ and 
$({\mathbf {\omega}}_i)_{i \not= j, j+1}$,   
\begin {eqnarray*}
\lefteqn {\sum_{{\mathbf {\omega}}_j, {\mathbf {\omega}}_{j+1}}
n_j r(z_j) \mu_{s_j, \lambda, \theta } ({\mathbf {\omega}}_j)
\mu_{s_{j+1}, \lambda, \theta } ({\mathbf {\omega}}_{j+1} ) 
f ( \Theta_{j+} ( {\mathbf {\omega}} )) 
}\\
&=&
\sum_{{\mathbf {\omega}}_j, {\mathbf {\omega}}_{j+1}}
n_{j+1} r(z_{j+1}) \mu_{s_j, \lambda, \theta } ({\mathbf {\omega}}_j)
\mu_{s_{j+1}, \lambda, \theta } ({\mathbf {\omega}}_{j+1})
f (  {\mathbf {\omega}} )
\end {eqnarray*}
and a similar dentity holds for the jumps to the left.
It follows that given the
local parities $n_j+z_j=s_j\mbox{ mod }2$, the probability
measures $\mu_{\us,\lambda,\theta}$ are stationary for the dynamics.
These are the \emph{equilibrium Gibbs measures} of our system. 
For a similar computation in the context of a simpler
one-component domain growth model see also \cite{balazs}.
% TO HERE 

Invariance under spatial translations is unfortunately lost
in this very general setup. In order to impose it, we restrict
ourselves to one of the following two choices: either $\us=\uzero$
or $\us=\uone$. 

\subsection{The hydrodynamic equations}

For the rest of this section we fix either
$\us=\uzero$ or $\us=\uone$ and we do not denote any more the
dependence on $\us$.

As we have mentioned already the globally conserved quantities
of our system are $\sum_j n_j$ and $\sum_j z_j$. In the
equilibrium regime $\mu_{\lambda,\theta}$ the averages of these
quantities are
% FROM HERE
%
\begin{equation*}
\rho
:=
\langle n_j \rangle_{\lambda,\theta}
=
\lambda\frac{\partial \log Z(\lambda,\theta)}{\partial \lambda},
\quad
\sigm
:=
\langle z_j \rangle_{\lambda,\theta}
=
\theta\frac{\partial \log Z(\lambda,\theta)}{\partial \theta}.
\end{equation*}
% TO HERE 
%
These are the particle density (per site) and the average slope
of the height of the wall, in equilibrium. 
It is easy to see that the map
$\R_+\times (1/\theta^*,\theta^*) 
\ni (\lambda,\theta) \mapsto (\rho,\sigm) \in
 \R_+\times\R$  
is globally invertible. Indeed, 
% FROM HERE
%
\begin{equation}
\left(
\begin{array}{cc}
{\partial \rho}/{\partial \lambda}
&
{\partial \rho}/{\partial \theta}
\\[3pt]
{\partial \sigm}/{\partial \lambda}
&
{\partial \sigm}/{\partial \theta}
\end{array}
\right)
=
\left(
\begin{array}{cc}
\var(n)
&
\cov(n,z)
\\[3pt]
\cov(n,z)
&
\var (z)
\end{array}
\right)
\left(
\begin{array}{cc}
\lambda^{-1}
&
0
\\[3pt]
0
&
\theta^{-1}
\end{array}
\right).
\label{invert}
\end{equation}
%
% TO HERE  
%
So the gradient matrix on the left hand side of (\ref{invert}) is
everywhere invertible and this implies global invertibility of
the map $(\lambda, \theta)\mapsto(\rho, \sigm)$.
With slight abuse of notation we denote the components of the
inverse function $\lambda=\lambda(\rho,\sigm)$ and
$\theta=\theta(\rho,\sigm)$.    From (\ref{symmetries}) it
follows that
\begin{equation}
\label{symmetries2}
\lambda(\rho,-\sigm)=\lambda(\rho,\sigm)
\quad\text{ and }\quad
\theta(\rho,-\sigm)=1/\theta(\rho,\sigm).
\end{equation}

In order to guess the system of hydrodynamic equations we have
to see first how the infinitisimal generator acts on the
conserved quantities. An easy computation shows:
%
% FROM HERE
%
\begin{align*}
&
Ln_j
=
\big( n_{j-1} r(z_{j-1}) - n_{j} r(- z_{j})  \big)
-
\big( n_{j} r( z_{j}) - n_{j+1} r(-z_{j+1})  \big)
\\
&
Lz_j
=
\big( n_{j-1} r(z_{j-1}) + n_{j} r(-z_{j}) \big)
-
\big( n_{j} r(z_{j}) + n_{j+1} r(-z_{j+1}) \big)
\end{align*}
%
% TO HERE
% 
On the right hand side of these equations we see \emph{discrete
gradients of fluxes}. This fact helps us guessing the
hydrodynamic equations.  Applying the standard formal
manipulations to our gradient system (see e.g. \cite{fritz},
\cite{kipnislandim}) and using the straightforward identities
%
% FROM HERE
%
\begin{equation*}
\langle n_j r(\pm z_j ) \rangle_{\lambda,\theta}
=
\lambda\theta^{\pm1}
\end{equation*}
in the hydrodynamic limit taken with \emph{hyperbolic (Eulerian)
scaling} of space and time, we arrive at the system of PDEs
\begin{equation}
\left\{
\begin{array}{l}
\pt \rho +
\px \big( \lambda(\rho,\sigm) 
(\theta(\rho,\sigm)-\theta(\rho,\sigm)^{-1}) 
    \big)
=0,
\\[5pt]
\pt \sigm +
\px \big( \lambda(\rho,\sigm) 
(\theta(\rho,\sigm)+\theta(\rho,\sigm)^{-1})
\big) 
=0.
\end{array}
\right.
\label{hydrodyneq}
\end{equation}

Under growth conditions on the rate function $r(z)$, as
$z\to\infty$, 
Yau's `relative entropy method' (see e.g. \cite{yau},  
%\cite{ollavaradhanyau}, 
\cite{fritz}, \cite{kipnislandim}) 
in principle can be applied to our system of interacting 
particles, resulting in the validity of the above system of PDEs in
the hydrodynamic limit, \emph{as long as the solutions are
smooth}. 
%
%
%  WORK IN PROGRESS 
%
% TO HERE

From the system (\ref{hydrodyneq}) we can derive the system 
(\ref{ourpde2})  by taking a second limit: 
We replace 
$\rho(t,x)$ by $\alpha^{2/3}\rho(\alpha t,\alpha^{2/3}x)$
and 
$\sigm(t,x)$ by $\alpha^{1/3} \sigm(\alpha t,\alpha^{2/3}x)$
We note  that  for small values of the variables $\rho$ and
$\sigm$, 
\begin{equation*}
\lambda(\rho, \sigm)=\rho +o(\rho),
\qquad
\theta(\rho,\sigm)=1+c \sigm + o(\sigm).
\end{equation*}
where 
\begin{equation*}
c=\left(\left.\frac{\partial^2 Z}{\partial
\theta^2}\right|_{\lambda=0,\theta=1} \right)^{-1}\in(0,\infty).
\end{equation*}
Letting now $\alpha\to0$, we arrive at (\ref{ourpde2}). 
We should emphasize here that this scaling limit does not depend
much on the details of microscopic system. Also, 
from any conservation law of the form 
\begin{equation*}
\left\{
\begin{array}{l}
\pt\rho + \px J(\rho,u) =0
\\[5pt]
\pt u   + \px K(\rho,u) =0
\end{array}
\right.
\end{equation*}
we would get (\ref{ourpde2}) under the same limiting procedure,
provided that 
\begin{equation*}
J(\rho,u)=\rho u +o(\rho u),
\quad
K(\rho,u)=\rho + o(\rho),
\quad
\text{ as }
\rho,u\to 0.
\end{equation*}
This indicates that (\ref{ourpde2}) is valid for a wider class 
of microscopic systems.

\section{Analysis of the PDE}

We are now going to see how the methods developed in the PDE
literature  (see \cite{hormander, serre, smoller}) can be applied
to our system. In order  to put things into perspective, we
briefly recall general 
results and see how they can be applied in the context of our
system (\ref{ourpde}).

\subsection{Two-component systems of hyperbolic conservation laws} 

For a generic two-component system we shall  use the notation
$v=v(t,x)=(v_1(t,x),v_2(t,x))^T$. (The superscript $^T$ will 
denote transposition of vectors/matrices.) The generic
two-component system is
\begin{equation}
\pt v + \px J (v) = 0,
\label{hcl}
\end{equation}
where $v\mapsto J(v)=(J_1(v),J_2(v))^T$ is a smooth vector field 
over $\R\times\R$. $J$ is the flux of the flow of 
the conserved quantity $v$. The initial conditions are specified by 
\begin{equation}
v(0,x)=v^{(0)}(x).
\label{incond}
\end{equation}

For a (possibly vector- or matrix valued)
function $f=f(v)$ we denote the gradient with respect to the
$v$-variables 
$\nabla f = 
({\partial f}/{\partial v_1},{\partial f}/{\partial v_2})$. 
For classical \emph{smooth} solutions $v(t,x)$, (\ref{hcl}) is 
equivalent to 
\begin{equation}
\pt v + \nabla J \cdot \px v = 0
\label{smoothhcl}
\end{equation}
(we use $\cdot$ to indicate products of matrices).

As a technical device one usually also considers the so-called
viscous equations
\begin{equation}
\pt v + \nabla J \cdot \px v = \vareps \px^2 v.
\label{visc}
\end{equation}
Existence and unicity of smooth solution $v^{(\vareps)}(t,x)$ of
(\ref{visc}), for any bounded and smooth initial conditions
(\ref{incond}) is
guaranteed by the smoothening effect of the \emph{artificial
viscosity term} on the right hand side. One hopes that physically 
acceptable (stable) solutions of the original system (\ref{hcl})
can be obtained as a \emph{strong} limit of the viscous solution
$v^{(\vareps)}(t,x)$, as $\vareps\to 0$. The existence of this
strong limit is a very difficult problem and is  
a main object of investigation in the context of hyperbolic
conservation laws.  

In our case (\ref{ourpde2}) the two components are
$
v_1=\rho$,
$v_2=u$,
and the corresponding fluxes are
$
J_1(\rho,u)=\rho u$,
$
J_2(\rho,u)=\rho
$. 
The inviscid system is (\ref{ourpde2}). The (artificially)
viscous system is
\begin{equation}
\left\{
\begin{array}{l}
\pt \rho+\px (\rho u) =\vareps \px^2 \rho
\\[5pt]
\pt u + \px \rho = \vareps \px^2 u.
\end{array}
\right.
\label{ourvisc}
\end{equation}
The viscous solutions (which do exist and are unique) will be
denoted by $\big(\rho^{(\vareps)}(t,x),
u^{(\vareps)}(t,x)\big)$.

\subsection{Hyperbolicity}

One has to check that the matrix $\nabla J$ has two distinct
real eigenvalues 
$\mu < \lambda$.
The domain where 
this 
holds will be denoted 
\begin{equation*}
{\cal D}_{\mathrm{hyp}}:=\{v\in\R\times\R: \mu(v)<\lambda(v)\}. 
\end{equation*}

The corresponding left (row) and right (column) eigenvectors
will be denoted by 
$l$ and $r$, respectively, $m$ and $s$. That is:
\begin{align}
&l \cdot\nabla J = \lambda l, 
\qquad 
&
\nabla J \cdot r = \lambda r,
\label{eveqlambda}
\\
&m \cdot \nabla J = \mu m,
\qquad
&
\nabla J \cdot s = \mu s.
\label{eveqmu}
\end{align}

For our system we find:
\begin{equation*}
\nabla J=
\left(
\begin{array}{cc}
u & \rho 
\\
1 & 0 
\end{array}
\right),
\end{equation*}
and
\begin{align}
&
\lambda=+\frac12(\sqrt{u^2+4\rho}+u),
\qquad
\phantom{m}
l=(\lambda,\rho),
\qquad
\phantom{s}
r=(\lambda,1)^T,
\label{ourevlambda}
\\
&
\mu=-\frac12(\sqrt{u^2+4\rho}-u),
\qquad
\phantom{l}
m=(\mu,\rho),
\qquad
\phantom{r}
s=(\mu,1)^T.
\label{ourevmu}
\end{align}
Note that 
$  l \cdot  s  =   m \cdot r = 0$,
as it should be. 

We conclude that for our system,
\begin{equation*}
{\cal D}_{\mathrm{hyp}}=
\{(\rho,u)\in\R\times\R:u^2+4\rho>0\}. 
\end{equation*}
Note that in the physically relevant 
domain with non-negative densities 
\begin{equation*}
{\cal D}_{\mathrm{ph}}:= 
\{(\rho,u)\in\R\times\R:\rho\ge0\},
\end{equation*}
there is one single point
where strict hyperbolicity is lost, namely $(\rho,u)=(0,0)$. 
On the other hand, we found that the system is still hyperbolic
in the physically meaningless domain 
${\cal D}_{\mathrm {hyp}} \setminus 
{\cal D}_{\mathrm {ph}}
= \{(\rho,u)\in{\cal D}_{\mathrm{hyp}}: \rho<0\}
\not=\emptyset$. 
At the moment nothing seems to prevent solutions to flow into 
this domain.
Later we
shall see that Lax's maximum principle (valid for stable entropy
solutions) takes care  of this problem. 

\subsection{Riemann invariants, characteristics}

In the generic two-component case, we are looking for scalar
functions 
${\cal D}_{\mathrm{hyp}}\ni v\mapsto w(v)\in \R$ 
and space-time trajectories
$\R_+\ni t\mapsto\xi(t)\in \R$ such that 
for smooth solutions of (\ref{hcl}) (or, equivalently, of 
(\ref{smoothhcl})) 
$w$ is conserved along the trajectory $\xi(t)$, i.e. 
\begin{equation*}
\frac{d}{dt} w\big( v(t,\xi(t) \big) = 0.
\end{equation*}
Using (\ref{smoothhcl}) we find:
\begin{equation}
\frac{d \xi}{d t}=
\frac{(\nabla w \cdot \nabla J)\cdot \px v}
     { \nabla w \cdot \px v }.
\label{wcons}
\end{equation}

In order to solve (\ref{wcons}), $\nabla w$ must be a left 
eigenvector of
the matrix $\nabla J$. It follows that this relation admits two 
solutions: one for each eigenvalue of $\nabla J$.
We denote the two solutions by $w$
(corresponding to the eigenvalue $\lambda$), respectively, by $z$
(corresponding to the eigenvalue $\mu$).  The gradients $\nabla
w$, respectively $\nabla z$, are parallel to the row vectors  $l$,
respectively $m$, defined in (\ref{eveqlambda}), respectively 
(\ref{eveqmu}). In other words,
% (see (\ref{leftdotright})), 
%
\begin{align*}
\nabla w \cdot s =0,
\quad
&
\frac{ d \xi }{d t} = \lambda,
\\
\nabla z \cdot r =0,
\quad
&
\frac{ d \xi }{d t} = \mu.
\end{align*}
These equations, of course, do not determine uniquely the
functions $w(v)$ and $z(v)$. Given two smooth, monotone maps
$f,g:\R\to\R$, the transformation $\hat w:=f(w)$, $\hat z:=g(z)$
leaves the above equations invariant. The functions $w$ and $z$
are called the {\it Riemann invariants}, or {\it characteristic
coordinates} of the problem. 

In our case the most convenient choice of the Riemann invariants
$w$ and $z$ is the following: let 
\begin{align*}
&
{\cal D}_w:=
\{ (\rho,u) \in {\cal D}_{\mathrm{hyp}} : \sqrt{u^2+4\rho}-u\ge0\}, 
\\
&
{\cal D}_z:=
\{ (\rho,u) \in {\cal D}_{\mathrm{hyp}} : \sqrt{u^2+4\rho}+u\ge0\},
\end{align*}
and define $w:{\cal D}_w\to\R$, $z:{\cal D}_z\to \R$ by the
formulas: 
\begin{align*}
&
w(\rho,u)
=
-\sqrt{\sqrt{u^2+4\rho}-u}
\left(\sqrt{u^2+4\rho}+2u\right),
\\
&
z(\rho,u)
=
-\sqrt{\sqrt{u^2+4\rho}+u}
\left(\sqrt{u^2+4\rho}-2u\right).
\end{align*}

\begin{figure}[ht]
\centering
%\vspace {-2cm}
\resizebox{0.95 \textwidth}{!} {\includegraphics {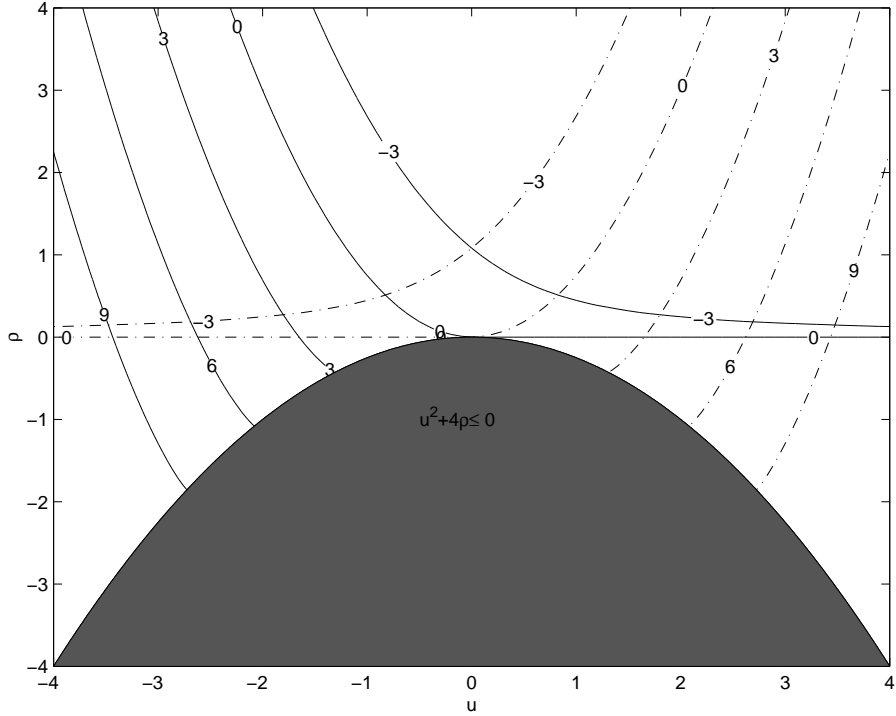}}
%\vspace {-1cm}
\caption{
Level lines of the Riemann invariants: $z=cst$ and $w=cst$}
\end{figure}

Note that 
${\cal D}_w\cap{\cal D}_z={\cal D}_{\mathrm{ph}}$, 
so that both  Riemann invariants are
defined in the physically relevant subdomain. 

It is straightforward to check that both Riemann invariants $w$
and $z$ defined above are \emph{convex} functions of the
variables $(\rho,u)$. This fact will have crucial importance in
later analysis. 

\subsection{Genuine nonlinearity}

In plain words, genuine nonlinearity of a two-component system
of hyperbolic conservation laws means that on the level curves
 $w(v)=\mbox{const.}$,
respectively $z(v)=\mbox{const.}$, the characteristic speed
$\mu$, respectively $\lambda$, varies strictly monotonically. 
Formally:
\begin{equation*}
\left. 
\frac{\partial \lambda}{\partial w} \right|_z
\not=0\not=
\left. \frac{\partial \mu}{\partial z} \right|_w.
\end{equation*}
Performing straightforward computations this turns out to be
equivalent to
\begin{equation*}
\nabla \lambda \cdot r 
\not=0\not=
\nabla \mu \cdot s.
\end{equation*}
That is: the characteristic speeds $\lambda$ and $\mu$ 
vary strictly monotonically in the direction of their corresponding  
right eigenvectors. 

In our case, given the formulas (\ref{ourevlambda}) and
(\ref{ourevmu}) we easily get
\begin{equation*}
\nabla \lambda \cdot r 
= \frac{2\lambda}{\lambda-\mu},
\qquad
\nabla \mu \cdot s 
= \frac{2\mu}{\mu-\lambda}.
\end{equation*}
Recall from (\ref{ourevlambda}), (\ref{ourevmu}) that 
on ${\cal D}_{\mathrm{ph}}$ we have $\mu\le0\le\lambda$,
with strict inequalities for $\rho>0$. 
We conclude that our system is genuinely nonlinear in the
interior of the physically relevant domain
${\cal D}_{\mathrm{ph}}$.
On the half lines $\rho=0$,
$u\le0$, respectively, $\rho=0$, $u\ge0$ (on the boundary 
of ${\cal D}_{\mathrm{ph}}$) genuine nonlinearity
of the first, respectively, of the second, characteristic speed
is lost.

\subsection{Weak solutions, shocks, Rankine-Hugoniot conditions}

As it is well-known, a nonlinear system of hyperbolic conservation laws 
(\ref{hcl}) can develop
singularities (e.g. discontinuities), irrespectively of the
smoothness of the initial conditions. A 
\emph{generalized} or \emph{weak} solution of (\ref{hcl}), (\ref{incond}) 
in a space-time domain is a
bounded, measurable function $(t,x)\mapsto v(t,x)$ satisfying
\begin{align}
\label{weaksln}
\int_{-\infty}^\infty\int_0^\infty
\big\{\pt\phi(t,x)
\cdot v(t,x)
+ 
\px\phi(t,x) \cdot
&
 J(v(t,x))\big\} dtdx
+
\\
&
+
\int_{-\infty}^\infty \phi(0,x)\cdot v^{(0)}(x) dx
=
0
\notag
\end{align}
for any row vector valued test function $\phi=(\phi_1,\phi_2)$
with compact support in the respective space-time domain. This
last equation is 
 obtained by a formal integration by parts. It is easily
seen that a strong (smooth) solution is also a weak solution.

Assuming a (locally) piecewise $C^1$ solution with a spatially
isolated jump discontinuity at some space-time position 
$(t,x)\in\R_+\times\R$, one derives the 
Rankine-Hugoniot conditions which relate the left- and right
limits of the function $x\mapsto v(t,x)$ at the discontinuity and the
propagation speed of the discontinuity:
\begin{equation}
\frac{J_1(v(t,x^+)) - J_1(v(t,x^-)) }{v_1(t,x^+) - v_1(t,x^-) }
=\s=
\frac{J_2(v(t,x^+)) - J_2(v(t,x^-)) }{v_2(t,x^+) - v_2(t,x^-) },
\label{rankine}
\end{equation}
where $\s$ is the propagation speed of the discontinuity, i.e.
the slope in space-time of the line of discontinuity. (\ref{rankine}) 
is derived from (\ref{weaksln}) by an elementary local 
argument, using the
divergence theorem (in space-time). Given the two independent
relations in (\ref{rankine}), any three of the five values 
$v_1(t,x^-)$, $v_1(t,x^+)$, $v_2(t,x^-)$, $v_2(t,x^+)$, $\s$ 
determine the
other two. This imposes a serious restriction on the possible
jump discontinuities of weak solutions. Note that the conditions
are left-right symmetric. 

We turn now to our system (\ref{ourpde2}). We denote by 
$(\rho^{\text{left}},u^{\text{left}})$,
respectively 
$(\rho^{\text{right}},u^{\text{right}})$, 
the values of the component functions
at the two sides of the presumed discontinuity. 
The Rankine-Hugoniot conditions are:
\begin{equation}
\label{ourrankine1}
\frac
{u^{\text{right}}\rho^{\text{right}}-u^{\text{left}}\rho^{\text{left}}}
{\rho^{\text{right}}-\rho^{\text{left}}}
=\s=
\frac
{\rho^{\text{right}}-\rho^{\text{left}}}
{u^{\text{right}}-u^{\text{left}}}.
\end{equation}
Given the value at one side of the discontinuity, the value at
the other side as function of propagation speed is expressed as
follows: 
\begin{equation}
\label{ourrankine2}
\rho^{\text{right}}=\s^2-\s u^{\text{left}},
\qquad
u^{\text{right}}=\s-\frac{\rho^{\text{left}}}{\s}
\end{equation}
Note that $\rho^{\text{right}}$, respectively,
$u^{\text{right}}$, 
is expressed as function of
$\s$ and $u^{\text{left}}$, respectively, as function of 
$\s$ and $\rho^{\text{left}}$,
only. (In principle, both should be expressed as functions of
$\s$, $\rho^{\text{left}}$ and $u^{\text{left}}$.) 
This is a special feature of our system.  

The propagation speed, as function of the values of the
components on both sides of the discontinuity, is expressed as:
\begin{equation*}
\s_{\pm}=
\pm\frac12
\big\{
\sqrt{(u^{\text{right}})^2+4\rho^{\text{left}}}\pm u^{\text{right}}
\big\}=
\mp\frac12
\big\{
\sqrt{(u^{\text{left}})^2+4\rho^{\text{right}}}\mp u^{\text{left}}
\big\}.
\end{equation*}

Lax's condition of stability for Rankine-Hugoniot  
discontinuities, \cite{lax1}, specified for two-component 
systems reads as follows: 
Assume that the weak solution (\ref{weaksln}) of the
two-component system (\ref{hcl}) is piecewise smooth, with a
spatially
isolated discontinuity with values $v^{\text{left}}$,
respectively,  $v^{\text{right}}$ on the two sides, propagating
according to the Rankine-Hugoniot conditions (\ref{rankine}). The
discontinuity is a stable \emph{back shock}, respectively,
\emph{front shock}, according whether
\begin{equation}
\label{backshock}
\mu(v^{\text{right}})
<\s<
\min\{\mu(v^{\text{left}}),\lambda(v^{\text{right}})\}, 
\end{equation}
or
\begin{equation}
\label{frontshock}
\max\{\lambda (v^{\text{right}}),\mu (v^{\text{left}})\}
<\s<
\lambda(v^{\text{left}}).
\end{equation}
Rankine-Hugoniot discontinuities which do not obey either one of
the conditions (\ref{backshock}) or (\ref{frontshock}), are
unstable, physically not realisable. 

Tedious (but, in principle straightforward) computations show,
that in the case of our system (\ref{ourpde}) the discontinuities
propagating according to (\ref{ourrankine1}), or equivalently
(\ref{ourrankine2}) are stable back shocks if $\s<0$ and stable
front shocks if $\s>0$.

\subsection{Entropies}

Given the two-component system of conservation laws (\ref{hcl}),
we look for \emph{additional} conserved quantities, i.e., for pairs
of functions 
${\cal D}_{\mathrm{hyp}}\ni v\mapsto (S(v),F(v))\in\R\times\R$ 
which satisfy
\begin{equation}
\pt S(v) + \px F(v)=0
\label{entropycons}
\end{equation}
for \emph{smooth solutions} of the original problem (\ref{hcl})
(or, equivalently: for smooth solutions of (\ref{smoothhcl})).
Indeed, (\ref{entropycons}) means, that $S(v(t,x))$ is globally 
conserved quantity, with flux $F(v(t,x))$. 
The pair of functions $(S,F)$
is called \emph{entropy/flux} pair. Using the form
(\ref{smoothhcl}), valid for smooth solutions of (\ref{hcl}), 
one finds the system of
PDEs defining an entropy/flux pair:
\begin{equation}
\nabla F = \nabla S\cdot \nabla J,
\label{entropyeq1}
\end{equation}
or, in extended form:
\begin{equation*}
\label{entropyeq2}
\frac{\partial F}{\partial v_k} 
=
\sum_{l=1}^2 
\frac{\partial S}{\partial v_l}
\frac{\partial J_l}{\partial v_k},
\qquad
k=1,2.
\end{equation*}
This is a two-component linear hyperbolic system of PDEs for the
two unknown functions $S$ and $F$ -- just well determined. There
are various alternative equivalent ways of writing it. E.g.,
eliminating the function $F$ we get a second order hyperbolic
PDE (a wave equation with variable  coefficients) for $S$:
\begin{equation*}
\label{entropyeq3}
\frac{\partial J_1}{\partial v_2}
\frac{\partial^2 S}{\partial v_1^2}
+
\left(
\frac{\partial J_2}{\partial v_2}-
\frac{\partial J_1}{\partial v_1}
\right)
\frac{\partial^2 S}{\partial v_1 \partial v_2}
+
\frac{\partial J_2}{\partial v_1}
\frac{\partial^2 S}{\partial v_2^2}
=0.
\end{equation*}
Or, changing variables  to the characteristic coordinates $(w,z)$:
\begin{equation*}
\label{entropyeq4}
\frac{\partial F}{\partial w} =
\lambda \frac{\partial S}{\partial z},
\qquad
\frac{\partial F}{\partial z} =
\nu \frac{\partial S}{\partial w}.
\end{equation*}
Or, eliminating $F$ between these two equations:
\begin{equation*}
\label{entropyeq5}
\frac{\partial^2 S}{\partial w\partial z}
=
\frac{1}{\lambda-\mu}
\left(
\frac{\partial \mu}{\partial w} 
\frac{\partial S}{\partial z}
-
\frac{\partial \lambda}{\partial z} 
\frac{\partial S}{\partial w}
\right).
\end{equation*}
These last two forms explicitly show the wave-character of the 
entropy equations (\ref{entropyeq1}).
Of particular importance are those entropy/flux pairs for
which the  function $v\mapsto S(v)$ is {\it convex}. 
 Such pairs will 
be simply called (with slight abuse of terminology) 
\emph{convex entropy/flux pairs}.

In the case of our system (\ref{ourpde2}) the entropy equations, 
written in terms of the physical variables $\rho$ and $u$,  are:
\begin{equation*}
\label{ourentropyeq1}
\frac{\partial F}{\partial \rho}=
u\frac{\partial S}{\partial \rho}+
\frac{\partial S}{\partial u},
\qquad
\frac{\partial F}{\partial u}=
\rho \frac{\partial S}{\partial \rho}
\end{equation*}
Or, eliminating $F$:
\begin{equation}
\label{ourentropyeq2}
\rho \frac{\partial^2 S}{\partial\rho^2}
-u \frac{\partial^2 S}{\partial\rho\partial u}
-\frac{\partial^2 S}{\partial u^2}
=0
\end{equation}

The existence of a  strictly convex  entropy/flux pair, 
\emph{globally defined} on 
${\cal D}_{\mathrm{ph}}=\{(\rho,u): \rho\ge 0, \ u\in\R\}$ and 
with $S$ bounded from below is very important, since the
applicability of Lax's Maximum Principle cited in the next
subsection relies on it. 
Here it is:
\begin{equation}
\label{ourfirstentropy}
S(\rho,u)=\rho\log\rho+\frac{u^2}{2},
\qquad
F(\rho,u)=u\rho(\log\rho+1).
\end{equation}

Lax's `entropy wave construction' (cf. \cite{lax2})
applies
also to our  system (\ref{ourpde}). Since these  
computations are rather involved, we do not reproduce them
here. Let us just point out, that this robust method ensures the
existence of a sufficiently rich family of convex entropy/flux
pairs in any fixed subdomain compactly contained in  ${\cal
D}_{\mathrm{ph}}$. 

There are also other (more ad hoc) methods of constructing 
entropy/flux pairs. Following, e.g.,  the ideas of 
\cite{lionsperthametadmor} we may try to find so called 
similarity solutions of the entropy equation
(\ref{ourentropyeq2}) 
of the form:
\begin{equation}
\label{similarity}
S(\rho,u)=\rho^{\alpha}\phi(\rho^\beta u)
\end{equation}
Elementary
computations show that $\beta=-1/2$ is the only choice
consistent with (\ref{ourentropyeq2}). 
Inserting (\ref{similarity}), with $\beta=-1/2$  
into (\ref{ourentropyeq2})
we find the following ordinary differential equation for the
function $\phi:\R\to\R$:
\begin{equation}
\label{similarityequ}
3(y^2-4/3)\phi''(y)+
(5-8\alpha)y\phi'(y)+
4\alpha(\alpha-1)\phi(y)=0.
\end{equation}
Any solution of (\ref{similarityequ}), with any $\alpha\in\R$
fixed provides an entropy of our system, via (\ref{similarity}).
So, we are able to construct a sufficiently rich family of
entropy/flux pairs to our system (\ref{ourpde2}).

\subsection{Entropy solutions}

A weak solution $(t,x)\mapsto v(t,x)$ of the generic system 
(\ref{hcl})
is called \emph{entropy solution} if for any convex entropy/flux pair
$(S,F)$ we have
\begin{equation}
\pt S +\px F \le 0
\label{entropycond}
\end{equation}
in the sense of distributions, i.e., for any \emph{positive} 
test function $(t,x)\mapsto \phi(t,x)$
\begin{equation*}
\int_{-\infty}^\infty\int_0^\infty
\left\{\pt \phi(t,x)S(v(t,x)) + \px \phi(t,x) F(v(t,x))\right\}
dtdx \ge 0.
\end{equation*}

Entropy solutions are the only physically admissible, stable
ones among the weak solutions. Strong limits of all convergent
approximation schemes (such as vanishing viscosity or various
convergent
finite difference schemes) result in entropy solutions. It is
also expected that convergent hydrodynamic limits of interacting
particle systems result in entropy solution of the corresponding
hyperbolic conservation laws. 
For piecewise smooth weak solutions, Lax's stability
condition for the shocks mentioned in a previous paragraph is
equivalent with the entropy conditions (\ref{entropycond}).

Of particular interest is the following Maximum Principle, due to 
P. Lax, see e.g. \cite{lax2}.

\smallskip

\noindent
{\bf Maximum Principle for Entropy Solutions.} 
\emph{
Assume that the following two conditions hold
\begin{enumerate}[(i)]
\item 
The Riemann invariants $v\mapsto w(v)$ and $v\mapsto z(v)$  
of the system of hyperbolic conservation laws (\ref{hcl}) 
are (globally) convex functions of $v$. 
\item
There exists  a globally defined convex entropy/flux pair, with
entropy function bounded from below.
\end{enumerate} 
Then, starting with bounded initial data,
$\sup_{-\infty<x<\infty}|v^{(0)}(x)|<\infty$, along entropy
solutions $(t,x)\mapsto v(t,x)$ the maximum values of the Riemann
invariants, 
$\sup_{-\infty<x<\infty}w(v(t,x))$
and 
$\sup_{-\infty<x<\infty}z(v(t,x))$ 
do not increase with $t$.
}

\smallskip

\noindent{\sl Remark:}
The same statement applies for solutions $v^{(\vareps)}(t,x)$ of
the viscous system (\ref{visc}) --- this follows from the
classical maximum principle. If $v^{(\vareps)}$ converges
strongly as $\vareps\to 0$, then the limiting $v$ is in fact an
entropy solution of the inviscid system (\ref{hcl}) and forcibly
it obeys Lax's Maximum Principle. It is not clear whether all
entropy solutions arise as limits of viscous solutions, with
vanishing viscosity. A general proof of the Maximum Principle for
entropy solutions can be found in \cite{lax2}.

Applying this theorem to our system we find that if we start with
bounded initial data $x\mapsto(\rho^{(0)}(x),u^{(0)}(x))\in{\cal
D}_{\mathrm{ph}}$ 
(that is: with non-negative initial density)  then \emph{entropy
solutions will stay in the physical domain, i.e., for any $t\ge0$
$x\mapsto(\rho(t,x),u(t,x))\in{\cal D}_{\mathrm{ph}}$.} 
(See Fig. 2 for graphical representation of the level curves
$w(\rho,u)=\text{const.}$ and $z(\rho,u)=\text{const}$.) This is
a very important consequence of the Maximum Principle: as we
already mentioned, a priori we could not see any reason banning a
(physically relevant) solution from flowing out into the
physically meaningless domain with $\rho<0$. 

In the case of isentropic gas dynamics, (\ref{igd}), choosing
convex versions of the Riemann invariants $w$ and $z$, for any
$w_{\text{max}}\in\R$, $z_{\text{max}}\in\R$, the domains 
\begin{equation*}
\{(\rho,m)\in\R_+\times\R: w(\rho,m)\le w_{\text{max}},
z(\rho,m)\le z_{\text{max}}\}
\end{equation*}
are compact. So starting with
bounded initial data global boundedness of (viscous and) entropy
solutions  is guaranteed by the  Maximum Principle. This is
unfortunately  not the case for our system. The domains
\begin{equation*}
\{(\rho,u)\in {\cal D}_{\mathrm{ph}}:
w(\rho,u)\le w_{\text{max}},
z(\rho,u)\le z_{\text{max}}\}
\end{equation*}
are \emph{not} compact, see Fig. 2. So here is an open question:
\emph{Is it the case, that if the initial data (\ref{ouric}) are
bounded then the solutions 
$\big(\rho^{(\vareps)}(t,x), u^{(\vareps)}(t,x)\big)$
of the viscous equation (\ref{ourvisc}) stay bounded for ever?
Similarly: is it the case that entropy solutions of
(\ref{ourpde2}) with bounded initial data stay bounded?} 
We guess that the answer to these questions are affirmative, but
we could  not prove this yet. 

\subsection{Vanishing viscosity, existence of entropy solutions}

The \emph{existence of entropy solutions} for a two-component
syetem of  hyperbolic conservationlaws (\ref{hcl}) is a
notoriously difficult question. The most powerful approach seems
to be the program initiated by R. DiPerna in \cite{diperna1}, 
completed for the case of isentropic gas dynamics (\ref{igd})
in \cite{diperna2}, then refined and extended in 
Lions et al. \cite{lionsperthametadmor} and in several other papers. 

In \cite{diperna1}, DiPerna proves the following result:

\smallskip

\noindent{\bf DiPerna's Theorem.}
\emph{
Consider the two-component system of hyperbolic
conservation laws (\ref{hcl}) and the corresponding viscous
system  (\ref{visc}). Assume that
\begin{enumerate}[(i)]
\item
The Riemann invariants $v\mapsto w(v)$ and $v\mapsto z(v)$ are
convex. (More precisely: there are convex choices of the Riemann
invariants. See subsection 4.3.)
\item
The system is genuinely nonlinear. (See subsection 4.4.)
\end{enumerate}
Let $\cal C$ be a domain compactly
contained in ${\cal D}_{\mathrm{hyp}}$ and assume that the
sequence of solutions $v^{(\vareps)}(t,x)$, $t\in[0,T]$,
$x\in\R$, $\vareps\to 0$,  of the viscous systems (\ref{visc}),
(\ref{incond}) 
takes values from $\cal C$. Then there is a subsequence
$v^{(\vareps^{'})}(t,x)$ which converges strongly in 
$L^1_{\mathrm{loc}}([0,T]\times\R)$. The limit $v(t,x)$ is
entropy solution of the system (\ref{hcl}). 
}

\smallskip

\noindent
{\sl Some Remarks:}
\begin{enumerate}[(1)]
\item
The proof relies on the construction of Lax's `entropy waves',
hinted at in subsection 4.6 and essentially on the so-called
compensated compactness method developed by Murat and Tartar. We
do not have a chance to reproduce here any technical part of the
proof. 
\item
It is assumed
that the viscous solutions stay in the domain $\cal
C$. 
However, even in this form the
theorem is technically very-very difficult. Extra difficulties
arise by relaxing this condition and imposing conditions {\em only on
the initial data}: in the isentropic gas dynamics and in our case
too, the solution data will typically flow to the boundary of the
domain of hyperbolicity and genuine nonlinearity, $\rho=0$, where
this theorem is not any more valid. 
\item
For extensions, physically more satisfactory formulations and
enormous further technical difficulties see e.g. \cite{diperna2},
\cite{lionsperthametadmor}, etc. 
\end{enumerate}

This theorem can be applied in a straightforward way for domains
$\cal C$, compactly contained in ${\cal D}_{\mathrm{ph}}$. We can
add to this that if initially 
\begin{equation*}
\max_x w(\rho^{(0)}(x), u^{(0)}(x))<0
\quad\mathrm{or}\quad
\max_x z(\rho^{(0)}(x), u^{(0)}(x))<0
\end{equation*}
then, due to the Maximum Principle, the viscous solutions 
$\rho^{(\vareps)}, u^{(\vareps)}$ are kept away from
the `dangerous' vacuum line $\rho=0$, see Fig. 2. So, in this
case one has to care only about the boundedness of the
solutions. 
 
\bigskip
\noindent
{\bf Acknowledegments.}
BT thanks illuminating discussions with 
M\'arton Bal\'azs, J\'ozsef Fritz and Benedek Valk\'o.  
We also thank Sophie Lemaire for kindly helping us producing Figure 2.
Cooperation between the authors is partially supported by the
French-Hungarian joint scientific research grant `Balaton'. 
   
%%%%%%%%%%%%%%%%%%%%%%%%%%%%%%%%%%%%%%%%%%%%%%%%%%%%%%%%%%%%%%%%%

-----------------------------------------------------------------
%%%%%%%%%%%%%%%%%%%%%%%%%%%%%%%%%%%%%%%%%%%%%%%%%%%%%%%%%%%%%%%%%

\vskip1cm

\hbox{\sc
\vbox{\noindent
\hsize66mm
B\'alint T\'oth\\
Institute of Mathematics\\
Technical University Budapest\\
Egry J\'oszef u. 1.\\
H-1111 Budapest, Hungary\\
{\tt balint@math.bme.hu}
}
\hskip5mm
\vbox{\noindent
\hsize66mm
Wendelin Werner\\
D\'ept. de Math\'ematiques\\
Universit\'e Paris-Sud\\
B\^at. 425\\ 
91405 Orsay cedex, France\\
{\tt wendelin.werner@math.u-psud.fr}
}
}

\end {document}